%% file: ffg-jap-revised.tex
\documentclass[12pt]{aptpub}
\usepackage{graphicx,subfigure}
\usepackage{amssymb}
\usepackage{alltt}

\usepackage{a4wide}
\usepackage{epsfig}
\usepackage{epsf}
\usepackage{float}
\usepackage{color}

\def\reff#1{(\ref{#1})}
\def\one{{\mathbf 1}}

\def\P{{\mathbb P}}
\def\R{{\mathbb R}}
\def\Z{{\mathbb Z}}
\def\square{\ifmmode\sqr\else{$\sqr$}\fi}
\def\sqr{\vcenter{
         \hrule height.1mm
         \hbox{\vrule width.1mm height2.2mm\kern2.18mm\vrule width.1mm}
         \hrule height.1mm}}                  

\def\uU{{\underline U}}
\def\ukappa{{\underline \kappa}}

\def\bx{{\mathbf x}}
\def\bz{{\mathbf z}}
\def\by{{\mathbf y}}
\def\bX{{\mathbf X}}
\def\bY{{\mathbf Y}}

\def\hG{{\hat G}}
\def\hX{{\hat X}}

\def\cB{{\mathcal B}}
\def\cF{{\mathcal F}}

\def\one{{\mathbf 1}}

\def\cF{{\mathcal F}}
\def\cN{{\mathcal N}}

\def\bX{{\mathbf X}}
\def\bY{{\mathbf Y}}
\def\ba{{\mathbf a}}
\def\bb{{\mathbf b}}
\def\bx{{\mathbf x}}
\def\by{{\mathbf y}}
\def\bSigma{{\bold \Sigma}}

\def\proof{{\noindent\bf Proof. }}

\authornames{P.J. Fern\'andez, P.A. Ferrari, S. P. Grynberg}
\shorttitle{Perfect sampling of multinormal distributions}

\begin{document}
\renewcommand{\baselinestretch}{1}
\title{Perfectly random sampling of truncated multinormal\\ distributions}

\authorone[Funda\c c\~ao Getulio Vargas]{Pedro J. Fern\'andez}

\authortwo[Universidade de S\~ao Paulo] {Pablo A. Ferrari}
\authorthree[Universidad de Buenos Aires]{Sebastian P. Grynberg}

\begin{abstract} 
The target measure $\mu$ is the distribution of a random vector in a box $\cB$,
a Cartesian product of bounded intervals. The Gibbs sampler is a Markov chain
with invariant measure $\mu$. A ``coupling from the past'' construction of the
Gibbs sampler is used to show ergodicity of the dynamics and to perfectly
simulate $\mu$. An algorithm to sample vectors with multinormal distribution
truncated to $\cB$ is then implemented. 
\end{abstract}


\section{Introduction}
\label{s1}

The use of latent variables has several interesting applications in
Statistics, Econometrics and related fields like quantitative marketing.
Models like tobit, and probit, ordered probit and multinomial probit are
good examples. References and examples could be found in Geweke, Keane and
Runkle (1997) and Allenby and Rossi (1999), especially for the multinomial
probit.

When the estimation process uses some simulation technique, in particular in
Bayesian analysis, the need of drawing a sample from the distribution of the
latent variable naturally arises. This procedure \emph{augments} the observed
data $Y$ with a new variable $X$ which will be referred as a \emph{latent}
variable. It is usually referred to as \emph{data augmentation} (Tanner and Wong
(1987), Tanner (1991), Gelfand and Smith (1990)) and in other contexts as
\emph{Imputation} (Rubin 1987).

The data augmentation algorithm is used when the likelihood function or
posterior distribution of the parameter given the latent data $\bX$ is simpler
than the posterior given the original observed data $\bY$.  If the distribution of
the variables $\bX$ is a multivariate normal in $d$ dimensions and $\bY=\one\{X\in
\cB\}$ is the indicator function of the set $\cB$, the drawing is made from the
normal distribution truncated to $\cB$ (and/or $\cB^c$).

The goal of this paper is to produce exact (or perfect) samples from random
variables with distributions supported on a $d$-dimensional box $\cB$; we call
\emph{box} the Cartesian product of $d$ bounded intervals.  We construct a
discrete-time stationary Markov process $(\zeta_t,\,t\in\Z)$ in the state space
${\cal B}$ whose time-marginal distribution at any time $t$ (that is, the law of
$\zeta_t$) has a given density distribution $g$ with support in~$\cB$.  The
construction is then implemented to perfectly simulate normal vectors of
reasonable large dimension truncated to bounded boxes. The approach is also
useful to show uniqueness of the invariant measure for a family of processes in
a infinite-dimensional space $[a,b]^{Z^d}$ with truncated Gaussian distributions
and nearest-neighbor interactions (in preparation, see also Section
\ref{s6}).

The process is the Gibbs sampler popularized by Geman and Geman (1985) and used
for the truncated normal case by Geweke (1991) and Robert (1995) in a
\emph{Markov Chain Monte Carlo (MCMC)} algorithm to obtain samples with
distribution approximating~$g$. To describe the process in our setting, let
$d\ge 2$ and take a box $\cB=B(1)\times\dots\times B(d)$ where $B(k)$ are
bounded intervals: $B(k) = [a_k,b_k]\subset\R$, $k=1,\dots,d$.  If
$\bx=(x_1,\dots,x_d)$ is the configuration at time $t-1$, then at time $t$ a
site $\kappa(t)$ is chosen uniformly in $\{1,\dots,d\}$, say $\kappa(t)=k$ and
$x_k$ is substituted by a value $y$ chosen with the density $g$ conditioned to
the values of the other coordinates; that is with the density $g_{k}(\cdot|\bx):
B(k)\to \R^+$ given by
\begin{equation}
  \label{p66}
g_{k}(y|\bx):= \frac{g({\by})}{\int_{B(k)} g({\bz})dz}\qquad y\in B(k)\,,
\end{equation}
where $k=\kappa(t)$, $\by=(x_1,\dots,x_{k-1},y,x_{k+1},\dots,x_d)$,
$\bz=(x_1,\dots,x_{k-1},z,x_{k+1},\dots,x_d)$. Note that $g_k(y|\bx)$ does not
depend on $x_k$. The distribution in ${\cal B}$ with
density $g$ is reversible measure for this dynamics.

We construct a stationary Gibbs sampler $(\zeta_t,\,t\in\Z)$ as a
function of a sequence $\uU=(U_t,t\in\Z)$ of independent identically distributed
and uniform in $[0,1]$ random variables and the updating schedule
$\ukappa=(\kappa(t),t\in\Z)$. We define (backwards) stopping times
$\tau(t)\in(-\infty,t]$ such that $\{\tau(t)=s\}$ is measurable with respect to
the field generated by the uniform random variables and the schedule in
$[s+1,t]$.  The configuration at time $t$ depends on the uniform random
variables and the schedule in $[\tau(t)+1,t]$.

In fact we construct simultaneously processes $(\zeta^\zeta_{[s,t]}, -\infty$
$<s\le t<\infty,\, \zeta\in\cB)$. For each fixed $s$ and $\zeta$,
$(\zeta^\zeta_{[s,t]},\,s\le t<\infty)$ is the Gibbs sampler starting at time
$s$ with configuration $\zeta$.  For each fixed $t$,
$(\zeta^\zeta_{[s,t]},\,s\le t,\,\zeta\in\cB)$ is a \emph{maximal coupling} (see
Thorisson (2000) and Section \ref{s2}). For each fixed time interval $[s,t]$ the
process $(\zeta^\zeta_{[s,t']}, s\le t'\le t)$ is a function of the uniform
random variables and the schedule in $[s+1,t]$ and the initial~$\zeta$. The
crucial property of the coupling is that $\zeta^\zeta_{\tau(t),t}$ does not
depend on $\zeta$.

The construction is a particular implementation of the \emph{Coupling From The
  Past (CFTP)} algorithm of Propp and Wilson (1996) to obtain samples of a law
$\nu$ as a function of a finite (but random) number of uniform random variables.
For an annotated bibliography on the subject see the web page of Wilson (1998).
A Markov process having $\nu$ as unique invariant measure is constructed as a
function of the uniform random variables. The processes starting at all possible
initial states run from negative time $s$ to time $0$ using the uniform random
variables $U_{s},\dots,U_0$ and a fixed updating schedule.  If at time 0 all the
realizations coincide, then this common state has distribution $\nu$.  If they
do not coincide then start the algorithm at time $s-1$ using the random
variables $U_{s-1},\dots,U_0$ and continue this way up to the moment all
realizations coincide at time 0; $\tau(0)$ is the maximal $s<0$ such that this
holds.  The difficulty is that unless the dynamics is monotone, one has to
effectively couple all (non countable) initial states.  Murdoch and Green (1998)
proposed various procedures to transform the infinite set of initial states in a
more tractable finite subset. 

In the normal case, for some intervals and
covariance matrices, the maximal coupling works without further treatment.
Moeller (1999) studies the continuous state spaces when the process undergoes
some kind of monotonicity; this is not the case here unless all correlations are
non negative.  Philippe and Robert (2003) propose an algorithm to perfectly
simulate normal vectors conditioned to be positive using coupling from the past
and slice sampling; the method is efficient in low dimensions.

Consider a $d$-dimensional random vector $\bY\in\R^d$ with normal distribution
of density
\begin{equation}
\label{To7}
f({\bx;\mu,\bSigma})=\frac{1}{(2\pi)^{d/2}|\bSigma|^{1/2}}e^{-Q({\bf
    x;\mu,\Sigma})/2}, \;\;\;\;{\bf x}\in \R^d, 
\end{equation}
where the covariance matrix ${\bSigma}$ is a positive definite matrix,
${\bf \mu}$ is the mean vector and
\begin{equation}
\label{To8}
Q({\bx;\mu,\bSigma})={(\bx -\mu)'\bSigma^{-1}(\bx-\mu)}.
\end{equation}
In this case we say that ${\bY}\sim {\cal N}({\mu,\bSigma})$.  We denote ${\bX}$
the truncation of ${\bY}$ to the box $\cal B$; its density is
\begin{eqnarray}
g({\bx};{\cB})=\frac{1}{Z({\cB})}f({\bx}){\bf 1}\{{\bx}\in {\cB}\},\label{p41}
\end{eqnarray}
where $Z({\cB})=\P[{\bX}\in {\cB}]=\int_{\cB}f({\bx})d{\bx}$ is the
normalizing constant.  

We have tested the algorithm in some examples. In particular for the matrix
$\bSigma^{-1} = \frac12 I +\frac12 {\mathbf 1}{\mathbf 1}'$ where $\mathbf I$ is the
identity matrix and $\mathbf 1' =(1,\dots,1)$. For the box $[0,1/2]^d$ and for
$d\ge 14$ our method is faster than the rejection method. The simulations
indicate that the computer time for the rejection method grows exponentially
with $d$ while it grows like $d^4$ for our method. Using a simple program in
Matlab in desktop microcomputer the method permits to simulate up to dimension
$d=30$ for the interval $[0,1/2]^d$. For bigger boxes or boxes apart from the
mean, for instance $[1/2,1]^d$, our method is less efficient but in any case it
is much faster than the rejection one. Since our method uses several times the
same uniform variables, to compare the computation time we have counted the
number of times each of the one-dimensional uniform random variable is used. It
is expected that if the box $\cB$ is sufficiently small and the covariance
matrix has a ``neighbor structure'', the number of steps may grow as $(d\ln
d)^2$.  In this case one should be able to simulate the state of an infinite
dimensional Gaussian field in a finite box.

Correlated normal vectors can be mapped into independent standard normal vectors
using the Cholesky transformation (see Chapter XI in Devroye (1986)). In the
case of truncated normals to a box, the transformation maps the box into a set.
One might simulate standard normals and reject if they do not belong to the
transformed box. The difficulties are similar to those of the rejection method:
as the dimension grows, the transformed set has small probability and the
expected number of iterations grows exponentially with the dimension.

When the truncating set is not a box, our method generally fails as the coupling
event for each coordinate has probability zero in general. We illustrate this
with an example in Section \ref{s6}.

Another possibility, also discussed by Devroye (1986) is to compute the marginal
of the first coordinate, then the second marginal conditioned to the first one
and so on. The problem here is that the computation of the marginals conditioned
to all the other coordinates be in $\cB$ may be as complicated as to compute the
whole truncated vector.

In Section \ref{s2} we define coupling and maximal coupling. In Section \ref{s3} we describe
the stationary theoretical construction of the Gibbs sampler and the properties
of the coupling. Section \ref{s4} is devoted to the pseudocode of the perfect
simulation algorithm. In Section \ref{s5} we compute the functions entering the
algorithm for the truncated normal case. In Section \ref{s6} we give some
examples and compare our perfect simulation algorithm with the rejection one
based on the uniform distribution.

\section{Coupling}
\label{s2}

A \emph{coupling} of a family of random variables $(X_\lambda,
\lambda\in \Lambda)$ with $\Lambda$ a label set, is a family of random
variables $(\hat{X}_\lambda: \lambda\in\Lambda)$ with the same
marginals; that is, such that
\[\hat{X}_\lambda\stackrel{D}{=}X_{\lambda}, \;\;\;\; \lambda\in\Lambda.\]
 An event $C$ is called \emph{coupling event}, if $\hat X_\lambda$ are
 identical in $C$. That is, 
 $$C\subseteq\{\hat{X}_\lambda=\hat{X}_{\lambda'} \; \mbox{for all } \;
 \lambda,\lambda' \in\Lambda\}.$$
We consider continuous random variables $X_\lambda$ in $\R$ with
densities $f_\lambda$ satisfying
\begin{eqnarray}
\label{Th2}
\P(C)\leq\int_{\R}\inf_{\lambda\in\Lambda}f_\lambda(y)\,dy.
\end{eqnarray}
for any coupling event $C$. This is always true if $\Lambda$ is countable and in
the normal case treated here.  When there is a coupling event $C$ such that
identity holds in \reff{Th2} the coupling is called \emph{maximal}. See
Thorisson (2000) for a complete treatment of coupling, including historical
quotations.

A natural maximal coupling of these variables is the following. Let $G_\lambda(x)
= \int_{-\infty}^x g_\lambda(y)dy$ be the corresponding cumulative distribution
functions and define
\begin{eqnarray}
\label{p20}
R(x):=\int_{-\infty}^{x}\inf_{\lambda\in \Lambda}g_\lambda(y)dy;\qquad R:=R(\infty)\,.
\end{eqnarray}
Let $\hG_\lambda(x) = G_\lambda(x)-R(x) + R$. Let $U$ be a random variable
uniformly distributed in $[0,1]$ and define
\begin{equation}
  \label{p23}
  \hX_\lambda =  
  R^{-1}(U)\one\{ U\le R\} + \hat G_\lambda^{-1}(U)\one\{ U> R\}\,,
\end{equation}
where $G^{-1}$, the \emph{generalized inverse} of $G$, is defined by
\begin{eqnarray}
G(x)\geq y \Leftrightarrow x\geq G^{-1}(y).
\end{eqnarray} 
The marginal distribution of $\hX_\lambda$ is $G_\lambda$:
\begin{eqnarray}
  \label{ppp}
  \P(\hX_\lambda<x) &=& \P(U<R(x) ) + \P(R<U<\hG_\lambda(x))
\;=\; R(x) + \hG_\lambda(x) - R 
\;=\; G_\lambda(x)\,.
\nonumber
\end{eqnarray}
The process $(\hX_\lambda,\lambda\in \Lambda)$ is a maximal coupling
for the family $(X_\lambda,\lambda\in \Lambda)$. We use this coupling
to construct the dynamics starting from different initial conditions.
Then we show how to compute $R(x)$ and $R$ in the normal case.

In the sequel we use a family of couplings. Let $(\Lambda_\ell,\,\ell\ge 1)$ be
a decreasing sequence of parameter sets: $\Lambda_{\ell+1}\subset\Lambda_{\ell}$
for all $\ell\ge 1$. Let $R_{[0]}(x)\equiv 0$ and for $\ell\ge 1$,
\begin{eqnarray}
\label{p20a}
R_{[\ell]}(x):=\int_{-\infty}^{x}\inf_{\lambda\in
  \Lambda_\ell}g_\lambda(y)dy;\qquad R_{[\ell]}:=R_{[\ell]}(\infty)\,. 
\end{eqnarray}
Assume $\lim_{\ell\to\infty}R_{[\ell]}= 1$. Define $D_{[\ell]}(x) =
R_{[\ell]}(x)-R_{[\ell-1]}(x)+ R_{[\ell-1]}$ and for $\lambda\in
\Lambda_\ell\setminus\Lambda_{\ell+1}$ define $\hat G_{\lambda,[\ell]}(x)=
G_\lambda(x) - R_{[\ell]}(x) + R_{[\ell]}$ and 
\begin{equation}
  \label{p23a}
  \hX_\lambda = \sum_{i= 1}^\ell D_{[i]}^{-1}(U) \,\one\{U\in[R_{i-1},R_i)\} \,+\,
  \hat G_{\lambda,[\ell]}^{-1}(U) \,\one\{U>R_{[\ell]}\} \,.
\end{equation}
Then $(\hX_\lambda,\,\lambda\in\Lambda_\ell)$ is a maximal coupling of
$(X_\lambda,\,\lambda\in\Lambda_\ell)$ for each $\ell$.

\section{A stationary construction of Gibbs sampler}
\label{s3}

The goal is to define a stationary process $(\zeta_t,\,t\in\Z)$ in $\cB$ with
the Gibbs sampler dynamics associated to a density $g$ in $\cB$. The time
marginal of the process at any time will be the distribution with the density
$g$. The Gibbs sampler is defined as follows.  Assume $\zeta_{t-1}\in\cB$ is
known, then at time $t$ choose a random site $\kappa(t)$ uniformly in
$\{1,\dots,d\}$ independently of ``everything else''. Then set $\zeta_{t}(k) =
\zeta_{t-1}(k)$ if $k\neq\kappa(t)$ and for $k=\kappa(t)$ choose $\zeta_t(k)$
with the law $g$ conditioned to the values of $\zeta_{t-1}$ at the other
coordinates:
  \begin{equation}
    \label{p60}
    \P(\zeta_{t} (k) \in A\,|\, \zeta_s, \,s< t)
=
\begin{array}{ll}
\int_A g_k(y|\zeta_{t-1}) dy &\hbox{if }\quad k= \kappa(t)\,,
\end{array}
\end{equation}
where $g_{k}(y|\bx)$, $y\in B(k):=[a_k,b_k]$, is the conditional density of the
$k$th component of $g$ given the other coordinates given in \reff{p66}.
Reversibility follows from the identity
\[
\int_\cB \int_{B(k)}H_1(\by)H_2(\bx)\, g_k(y|\bx)\,dy\,g(\bx) d\bx =
\int_{\cB}\int_{B(k)} H_1(\bx)H_2(\by)\, g_k(y|\bx)\, dy\, g(\bx)d\bx
\]
for any $H_1$ and $H_2$ continuous functions from $\cB$ to $\R$, where
$\bx=(x_1,\dots,x_d)$ and $\by=(x_1,\dots,x_{k-1},y,x_{k+1},\dots,x_d)$.  In
particular the law with density $g$ is invariant for the Gibbs sampler. It is
indeed the unique invariant measure for the dynamics, as we show in Theorem
\ref{p31}.

To define a stationary version of the process we construct a family of couplings
for each time interval $[s,t]$. The coupling starts with all possible
configurations at time $s$ and uses random variables associated to each time to
decide the updating. All couplings use the same updating variables.  To reflect
the evolution of all possible configurations we introduce a process on the set
of boxes contained in $\cB$. The initial configuration of this process is $\cB$
at time $s$ but at times $s'>s$ can be (and will be as $s$ decreases) reduced
to a point.  Let $\xi$ be a $d$-tuple of closed bounded intervals of $\R$,
$\xi(i)\subset B(i)$, $i=1,\dots,d$. We abuse notation calling $\xi$ also the
box
\begin{eqnarray}
\nonumber
\xi=\xi(1)\times\cdots\times\xi(d)\subset \cB.
\end{eqnarray} 
Notice that $\xi$ has dimension less than $d$ when $\xi(i)$ is just a point for
some $i$.
Define
\begin{eqnarray}
\label{integrales1}
R_k(x|\xi)&=&
\int_{a_k}^{x}
\big(\min_{\bx\in\xi}g_{k}(y|\bx)\big)\,dy,\qquad x\in B(k)\\ 
\label{integrales2}
R_k(\xi)&=&R_k(b_k|\xi)\,.
\end{eqnarray}
These functions depend on $\xi$ only through $(\xi_i, \,i\neq k)$. $R_k(\xi)$ is
the probability that if the $k$-th coordinate is updated when the other
coordinates belong to the set $\xi$, then a coupled event is attained for the
$k$ coordinate. This means that for any configuration $\bx\in\xi$, the updated
value of the $k$-th coordinate will be the same, say $x$; its law is given by
$R_k(x|\xi)$. In this case the updated interval $\xi(k)$ is also reduced to the
point $x$. If the coupled event is not attained $\xi(k)$ is kept as the interval
$B(k)$.

The updating random variables consist on two families: $\uU=(U_t : t\in\Z)$, a family of independent variables with uniform distribution in $[0,1]$ and
$(\kappa(t) :\, t\in\Z)$, a family of independent variables with uniform
distribution in $\{1,\dots,d\}$ and independent of $\uU$. 

Now for each $s\in\Z$ we define a process $(\eta_{[s,t]},\,t\ge s)$ taking
values on boxes contained in $\cB$ as function of $((U_t,\kappa(t)),\, t>s)$.
The initial state is $\eta_{[s,s]}=\cB$ and later each coordinate $k$ is either
a (random) point or the full interval $B(k)$. More precisely, for $s\in\Z$ and
$x\in \R$ set
\begin{eqnarray}
R_{[s,s]}= 0\,, 
\qquad
R_{[s,s]}(x)= 0\,, 
\qquad D_{[s,s]}(x)=0\,,
\qquad
\eta_{[s,s]}=\cB.
\end{eqnarray}
Fix $n\ge 1$ and assume $R_{[s,t]}(x)$, $R_{[s,t]}$ and $\eta_{[s,t]}$
are defined if $0\le t-s\le n-1$. Then for $t-s= n$ and $x\in
B(\kappa(t))$, define
\begin{eqnarray}
R_{[s,t]}&=&R_{\kappa(t)}(\eta_{[s,t-1]}) \\
R_{[s,t]}(x)&=&R_{\kappa(t)}(x|\eta_{[s,t-1]}), \\
D_{[s,t]}(x)&=&R_{[s+1,t]}+R_{[s,t]}(x)-R_{[s+1,t]}(x)\,\nonumber
\\\label{oi}
\eta_{[s,t]}(k)&=&\left\{
\begin{array}{lll}
\eta_{[s,t-1]}(k)\,,\qquad \mbox{if}\;k\neq \kappa(t),\\
\\
\eta_{[s+1,t]}(k) \,{\bf 1}\{U_t\leq
R_{[s+1,t]}\}\,+ \,D_{[s,t]}^{-1}(U_t)\,{\bf 1}\{R_{[s+1,t]}<U_t\leq
R_{[s,t]}\}\,
\\ \qquad\qquad\qquad +\,B(k)\,{\bf 1}\{U_t> R_{[s,t]}\}\,,
\qquad\mbox {if}\;k=\kappa(t);
\end{array}
\right.\,.
\end{eqnarray}
In words: $k=\kappa(t)$ defines the coordinate to be updated at time $t$.
$R_{[s,t]}$ is the probability that the coupling event is attained at coordinate
$\kappa(t)$ for all the processes starting at times smaller than or equal to
$s$. The coupling event is attained for the process starting at $s$ when
$U_t<R_{[s,t]}$. In case the coupling event is attained, the value of coordinate
$k$ is given by $D^{-1}_{[s',t]}(U_t)$, for $s'$ given by the biggest
$s\le t$ such that $U_t<R_{[s,t]}$ (second term in \reff{oi}).  This value is
the same for each $s\le s'$ (first term in \reff{oi}).  When the coupling
event is not attained (that is, for $s>s'$), the coordinate $k$ is kept
equal to the full interval $B(k)$ (third term in \reff{oi}).

It follows from this construction that
$(R_{[t-1,t]}(\cdot),\,R_{[t-1,t]},\,D_{[t-1,t]})$ does not depend on $\uU$ and,
for $s<t-1$,
\begin{eqnarray}
  \label{p64}
  &&( R_{[s,t]}(\cdot),\,R_{[s,t]},\,D_{[s,t]})
  \hbox{ is a function of }
  ((U_n,\kappa(n)),\,n=s+1,\dots,t-1).\\
  &&\eta_{[s,t]} \hbox{ is a function of }
  (R_{[s,t]}(\cdot),\,R_{[s,t]},\,D_{[s,t]}) \hbox{ and } (U_t,\kappa(t))\,.
\end{eqnarray}
 The $k$th coordinate of
$\eta_{[s,t]}$ is either the interval $B(k)$ or a point. The process
is monotonous in the following sense:
\begin{equation}
  \label{p30}
  \eta_{[s-1,t]}\subset\eta_{[s,t]}\hbox{ for all } s<t\,.
\end{equation}
In particular, if some coordinate is a point at time $t$ for the
process starting at $s$, then it will be the \emph{same} point for the
processes starting at previous times:
\begin{equation}
  \label{p29}
 \hbox{If }\eta_{[s,t]}(k) \hbox{ is a point, then }\eta_{[s-1,t]}(k)
=\eta_{[s,t]}(k),\qquad \hbox{ for all }s<t\,.
\end{equation}
For each time $t$ define
\begin{equation}
  \label{p29a}
  \tau(t) = \max\{s<t\,:\, \eta_{[s,t]}(k) \hbox{ is a point for all }k\}\,. 
\end{equation}

Using \reff{p64} we conclude that $\tau(t)+1$ is a stopping time for
the filtration generated by $((U_{t-n},\kappa(t-n)),n\ge 0)$: the
event $\{\tau(t) = s\}$ is a function of $((U_n,\kappa(n)),\,
n=s+1,\dots,t)$. 

Assume $\tau(t)>-\infty$ almost surely for all $t$ and define the process
$(\zeta_t,\, t\in\Z)$ in $\cB$ by
\begin{equation}
  \label{p30b}
  \zeta_t := \eta_{[\tau(t),t]} \in \cB\,.
\end{equation}
Noticing that $\tau(t)>-\infty$ is equivalent to $R_{[n,t]}\to 1$ as
$n\to-\infty$, we get the following explicit expression for $\zeta_t$:
\begin{eqnarray}
\label{p68}
\zeta_t(k)&=&\left\{
\begin{array}{lll}
\zeta_{t-1}(k)\,,\qquad \mbox{if}\;k\neq \kappa(t),\\
\\
\sum_{n\le t-1} \,D_{[n,t]}^{-1}(U_t)\,{\bf 1}\{R_{[n+1,t]}<U_t\leq
R_{[n,t]}\}\,,
\qquad\mbox {if}\;k=\kappa(t); 
\end{array}
\right.\,.
\end{eqnarray}

For each $s<t$ and $\zeta\in\cB$ we now construct the process
$(\zeta^\zeta_{[s,t]},\, t\ge s)$, the Gibbs sampler starting with $\zeta$ at
time $s$.  Let $G_k(x|\bx)=\int_{a_k}^x g_k(z|\bx)dz$ and
\begin{equation}
  \label{p33}
  \hat G_{[s,t]}(x|\zeta) = R_{[s,t]}+ G_{\kappa(t)}(x|\zeta) - R_{[s,t]}(x)\,.
\end{equation}
For each $\zeta\in\cB$ define $\zeta^\zeta_{[s,s]}=\zeta$ and for $t>s$,
\begin{eqnarray}
\label{p67}
\zeta^\zeta_{[s,t]} (k)&=&\left\{
\begin{array}{lll}
  \zeta^\zeta_{[s,t-1]}(k)\,,\qquad \mbox{if}\;k\neq \kappa(t),\\
  \\
  \sum_{n=s}^{t-1}
  \one\{U_t\in[R_{[n+1,t]},R_{[n,t]}]\}\, D^{-1}_{[n+1,t]}(U_t)\\
\qquad\qquad
  \qquad\qquad +\;
  \one\{U_t>R_{[s,t]}\}
  \,\hG_{[s,t]}^{-1}(U_t|\zeta^\zeta_{[s,t-1]})\,,
  \qquad\mbox { if }\;k=\kappa(t)
\end{array}
\right.\,.
\end{eqnarray}
For each fixed $t$ and $k$, $((\zeta^\zeta_{[s,t]}(k),\, \zeta\in\cB),\,s\le t)$
is a maximal coupling among the processes starting with all possible
configurations $\zeta$ at all times $s\le t$.

\begin{theorem}
  \label{p31}
  Assume $\tau(t)>-\infty$ almost surely for all $t\in\Z$. Then the process
  $(\zeta_t,\,t\in\Z)$ defined in \reff{p68} is a stationary Gibbs sampler
  related to $g$.  The distribution of $\zeta_t$ at each time $t$ is the
  distribution with density $g$ in $\cB$. This distribution is the unique
  invariant measure for the process and
  \begin{equation}
    \label{q77}
    \sup_{\zeta} \P(\zeta^\zeta_{[s,t]} \neq \zeta_t) \le \P(\tau(t)<s)\,,
  \end{equation}
  where $\zeta^\zeta_{[s,t]}$ is the process starting with configuration $\zeta$
  at time $s$.
\end{theorem}

\proof Stationarity follows from the construction. Indeed, for all
$s\le t$ and $\ell\in\Z$, $\eta_{[s,t]}$ as a function of
$((U_n,\kappa(n)), n\in \Z)$ is identical to $\eta_{[s+\ell,t+\ell]}$
as a function of $((U_{n+\ell},\kappa(n+\ell)), n\in\Z)$.

The function $\Psi_t:[0,1]\times [0,1]^{\{-\infty,\dots,t-1\}}\to
B(\kappa(t))$ defined in \reff{p68} by
\begin{equation}
  \label{p69}
  \Psi_t(U;U_{t-1},U_{t-2}, \dots) =\sum_{n\le t-1} \,D_{[n,t]}^{-1}(U)\,{\bf
    1}\{R_{[n+1,t]}<U\leq R_{[n,t]}\} 
\end{equation}
satisfies for $k=\kappa(t)$
\begin{equation}
  \label{p70}
  \P(\Psi_t(U;U_{t-1},U_{t-2}, \dots)\in A) = \int_A g_{k} (y|\zeta_{t-1})dy\,,
\end{equation}
where we used that $\zeta_{t-1}$ is function of $U_{t-1},U_{t-2},\dots$ and
$\kappa_{t-1},\kappa_{t-2},\dots$. This is sufficient to show that $\zeta_t$
evolves according to the Gibbs sampler.

To finish we need to show that the distribution with density $g$ is the marginal
law of $\zeta_t$. Since the updating is performed with the conditioned
distribution, the distribution with density $g$ is invariant for the Gibbs
sampler.

Since $\zeta^\zeta_{[s,t]}=\zeta_{[s,t]}$ if $s\le \tau(t)$, then \reff{q77}
follows. This also implies that $\lim_{s\to-\infty}\zeta^\zeta_{[s,t]} =\zeta_t$
almost surely. In particular $\zeta^\zeta_{[s,t]}$ converges weakly to the law
of $\zeta_t$ as $s\to-\infty$.  The same is true for
$\lim_{t\to\infty}\zeta^\zeta_{[s,t]}$. Taking $\zeta$ random with law $g$, this
proves $\zeta_t$ must have law $g$. Then taking $\zeta$ with law $\tilde g$
invariant for the dynamics, we conclude $\tilde g =g$ and the uniqueness of the
invariant measure follows.  \square

\begin{lemma}
  \label{s22}
A sufficient condition for $\tau(t)>-\infty$ \emph{a.s.} is that $R_k(\cB)>0$
for all $k$. In this case $\P(\tau(t)<s)$ decays exponentially fast with $t-s$. 
\end{lemma}

\proof $\tau(t)>\tau^0(t) := \max\{s<t-d\,:\,$ $\kappa(s+n) = n$ and
$U_{s+n}<R_{\kappa(s+n)}(\eta), n=1,\dots,d\}$, the last time in the past the
coordinates have been successively updated to a point independently of the other
coordinates. At $\tau^0(t)+d$ and successive times $t'\ge\tau^0(t)+d$, each
coordinate of $\eta_{[\tau^0(t),t']}$ has been reduced to a point. $\tau_0(t)$
is finite for all $t$ because the event $\{\kappa(s+n) = n$ and
$U_{s+n}<R_{\kappa(s+n)}(\eta), n=1,\dots,d\}$ occurs for infinitely many $s$
with probability one. The same argument shows the exponential decay of the tail
of $\tau(t)$.  \square

\paragraph{Remark} The velocity of convergence of the Gibbs sampler to
equilibrium depends on the values $(R_k(\cB),\,k=1,\dots,d)$. In turn, these
values depend on the size of the box $\cB$ and the correlations of the
distribution $g$. Strongly correlated vectors produce small values $R_k(\cB)$
and hence slow convergence to equilibrium.  The efficiency of the algorithms
discussed in the next sections will also depend on these values.

\section{Perfect simulation algorithm}
\label{s4}
The construction of Section \ref{s3} is implemented in a perfect simulation
algorithm.  Let $\xi$ be a box contained in $\cB$. For each $k=1,\dots,d$, let
$R_{k}(\cdot|\xi)$ and $R_{k}(\xi)=R_{k}(b_{k}|\xi)$ as in \reff{integrales1}
and \reff{integrales2}.  When $\xi=\cB$ we denote by $R_{k}$ the value of
$R_{k}(\cB)$.  Let $\kappa(t)$ be the updating schedule. It can be either a
family of iid chosen uniformly in $\{1,\dots,d\}$ or the periodic sequence
$\kappa(t)= [t-1]_{{\rm mod }\ d}+1$.  For each pair $\eta,\xi$ of boxes
contained in $\cB$ and $1\leq k\leq d$, let $D_{k,\eta,\xi}:\cB(k)\to[0,1]$ be the function
defined by
\begin{eqnarray*}
D_{k,\eta,\xi}(x)=R_{k}(x|\xi)-R_{k}(x|\eta)+R_{k}(\eta)\;\;\;\;\;(x\in\cB(k))\,.
\end{eqnarray*}

Let $\xi_1,\xi_2$ be boxes contained in $\cB$, $u\in[0,1]$,
and $1\leq k\leq d$. The auxiliary {\em coupler function} $\phi$ is defined by
\begin{alltt}
function \(\phi(\xi\sb{1},\xi\sb{2},u,k)\):
\(\xi\leftarrow\xi\sb{1}\)
\(\xi(k)\leftarrow\xi\sb{2}(k)\)
if \(\xi\sb{2}(k)\) is not a point in \(\R\) and \(u\leq R\sb{k}(\xi)\),
   update the \(k\)th coordinate of \(\xi\):  
   \(\xi(k)\leftarrow D\sb{k,\xi\sb{2},\xi}\sp{-1}(u)\)
end
\(\xi\sb{2}\leftarrow\xi\)
Return(\(\xi\sb{2}\))
\end{alltt}

\paragraph{\bf Perfect Simulation Algorithm:} Perform the following steps

\begin{alltt}
\(T\leftarrow 0\)
start(\(\eta\sb{0}\)) 
while \(\eta\sb{0}\) is not a point in \(\cB\)
      \(T\leftarrow T-1\)
      start(\(\eta\sb{T}\))
      \(t\leftarrow T\)
      while \(t< 0\) update \(\eta\sb{t+1}\) using the coupler function \(\phi\):
            \(\eta\sb{t+1}\leftarrow\phi(\eta\sb{t},\eta\sb{t+1},U\sb{t+1},\kappa(t+1))\)
            \(t\leftarrow t+1\)
      end      
end
Return(\(\eta\sb{0}\))
\end{alltt}

Where, for each $T\leq 0$, ``start($\eta_T$)'' is the following subroutine

\begin{alltt}
start(\(\eta\sb{T}\)):
\(\eta\sb{T}\leftarrow\cB\)
Generate \(U\sb{T}\) a uniform random variable in \([0,1]\)
if \(U\sb{T}\leq R\sb{\kappa(T)}\), compute the \(\kappa(T)\)th coordinate of \(\eta\sb{T}\):
   \(\eta\sb{T}(\kappa(T))\leftarrow R\sb{\kappa(T)}\sp{-1}(U\sb{T}|\eta\sb{T})\)
end
\end{alltt}

\section{Truncated Normal distributions}
\label{s5}

In this section we implement the construction in the normal case. We start with
elementary facts of the one dimensional Normal distribution.

\paragraph{One dimension}
Let $\phi(x)=\frac{1}{\sqrt{2\pi}}e^{-x^2/2}$ be the standard normal
density in $\R$ and $\Phi(x)=\int_{-\infty}^x\phi(y)dy$ the
corresponding distribution function. Let $a<b$ be real numbers.
A random variable $X$ with mean $\mu$ and variance $\sigma^2$ has
\emph{truncated normal distribution} to the interval $[a,b]$ if its
density is given by
\begin{equation}
\label{S6}
g_{a,b}(x;\mu,\sigma^2)=\frac{1}{\sigma}\frac{\phi(\frac{x-\mu}{\sigma})}
{\Phi(\frac{b-\mu}{\sigma})-\Phi(\frac{a-\mu}{\sigma})},
\;\;\;\; a\leq x\leq b, 
\end{equation}
where $\mu\in\R$ and $\sigma>0$. This distribution is called
${\cal T}_{a,b}{\cal N}(\mu,\sigma^2)$. The truncated normal is the
law ${\cal N}(\mu,\sigma^2)$ conditioned to $[a,b]$. 
Let
\begin{eqnarray}
\label{S1}
A_{a,b}(\mu)=\Phi(b-\mu)-\Phi(a-\mu).
\end{eqnarray}
Let $X_{\mu}, \mu\in I=[\mu^-,\mu^+]$ be a family of normal distributions
truncated to the interval $[a,b]$: $X_{\mu}\sim {\cal T}_{a,b}{\cal
  N}(\mu,\sigma^2)$.  Let $A^\pm=A^\pm_{a,b}$ and
$x(I,\sigma)=x_{a,b}(I,\sigma)$ be defined by
\begin{eqnarray}
\label{S13b}
A^-=\Phi\Bigl(\frac{b-\mu^-}{\sigma}\Bigr)-\Phi\Bigl(\frac{a-\mu^-}{\sigma}\Bigr),
\;\;\;\;
A^+=\Phi\Bigl(\frac{b-\mu^+}{\sigma}\Bigr)-\Phi\Bigl(\frac{a-\mu^+}{\sigma}\Bigr)  
\end{eqnarray}
and
\begin{equation}
\label{S13}
x(I,\sigma)=\frac{\mu^-+\mu^+}{2}-\frac{\sigma^2}{\mu^+-\mu^-}\ln
\left(\frac{A^-}{A^+}\right)\,. 
\end{equation}
Let $I=[\mu^-,\mu^+]$ and 
\begin{eqnarray}
  \label{p40}
  R(x|I) &=& \int_{a}^x\inf_{\mu\in I} g_{a,b}(y;\mu,\sigma^2)dy, \qquad x\in[a,b]\\
R(I) &=& R(b|I)\,.
\end{eqnarray}
The next proposition, proven later, gives an explicit way of computing the
functions participating in the multicoupling. It says that the infimum of normal
densities with same variances $\sigma$ and means in the interval $I=[\mu^-,\mu^+]$
coincides with the normal with mean $\mu^-$ up to $x(I,\sigma)$ and with the normal
with mean $\mu^+$ from this point on.
\begin{proposition}
\label{p22}
The integrand in \reff{p40} has the following explicit expression
\begin{equation}
\label{Marzo1}
\inf_{\mu\in I}g_{a,b}(x,\mu,\sigma^2)=g_{a,b}(x;\mu^-,\sigma^2)
\one\{x(I,\sigma)\leq x\} 
+g_{a,b}(x;\mu^+,\sigma^2) \one\{x(I,\sigma)>x\}.
\end{equation}
As a consequence,
\begin{eqnarray}\label{S44}
R(x|I)
&=&\frac{1}{A^{+}}\left[\Phi\left(\frac{x-\mu^+}{\sigma}\right)
-\Phi\left(\frac{a-\mu^+}{\sigma}\right)\right]{\bf 1}\{x<x(I,\sigma)\}\nonumber\\
&&\qquad+\,\frac{1}{A^{+}}\left[\Phi\left(\frac{x(I,\sigma)-\mu^+}{\sigma}\right)
-\Phi\left(\frac{a-\mu^+}{\sigma}\right)\right]{\bf 1}\{x\geq x(I,\sigma)\}\nonumber\\
&&\qquad+\,\frac{1}{A^{-}}\left[\Phi\left(\frac{x-\mu^-}{\sigma}\right)
-\Phi\left(\frac{x(I,\sigma)-\mu^-}{\sigma}\right)\right]{\bf
1}\{x\geq x(I,\sigma)\}.
\nonumber
\end{eqnarray}
In particular, 
\[
R(I)\;=\;\frac{1}{A^{+}}\left[\Phi\left(\frac{x(I,\sigma)-\mu^+}{\sigma}\right)
  -\Phi\left(\frac{a-\mu^+}{\sigma}\right)\right]+
\frac{1}{A^{-}}\left[\Phi\left(\frac{b-\mu^-}{\sigma}\right)
  -\Phi\left(\frac{x(I,\sigma)-\mu^-}{\sigma}\right)\right]
\] 
and for each $u\in [0,R(I)]$ there exists a unique $x_u\in [a,b]$ such
that $R(x_u|I)=u$.
\end{proposition}

\paragraph{Remark}  $0<R(I)<1$. Indeed, $R(I)>0$ because 
$\inf_{\mu\in I} g_{a,b}(x;\mu,\sigma^2)>0$. 
If $R(I)=1$, then all densities coincide unless
$\mu^-=\mu^+$; by hypothesis this trivial case is excluded.

\paragraph{Multivariate Normal}

Let ${\bf Y} \sim{\cal N}({\bf \mu,\bSigma})$, ${\bf X}= {\bf Y|_{\cal B}}$ be
the vector conditioned to ${\cal B}$ and ${\bf x}\in {\cal B}$. The law of $X_k$
conditioned to $(X_i,\, i\neq k)=(x_i,\,i\neq k)$ is the truncated normal ${\cal
  T}_{a_k,b_k}{\cal N} (\hat{\mu}_k({\bf x}),\hat{\sigma}_k^2)$,
where 
\begin{eqnarray}
\label{I2}
\hat{\mu}_k({\bx})=\mu_k-\frac{1}{a_{k k}}\sum_{i\neq k}a_{k i}(x_i-\mu_i),
\;\;\;\; \hat{\sigma}_k^2=\frac{1}{a_{k k}}. 
\end{eqnarray}
We see that $\hat \mu_k(\bx)$ does not depend on $x_k$ and 
\begin{eqnarray*}
  \hat{\mu}_k({\bf x})\in [\mu_k^-({\bf x}),\mu_k^+({\bf x})]\,.
\end{eqnarray*}
where
\begin{eqnarray}\nonumber
\mu_k^-({\bf x})=\frac{1}{a_{k k}}\sum_{i=1}^da_{k
  i}\mu_i-\frac{1}{a_{k k}}\left(\sum_{a_{k i}<0, i\neq k}a_{k
    i}a_i+\sum_{a_{k i}\geq 0, i\neq k}a_{k i}b_i\right),\\ 
\mu_k^+({\bf x})=\frac{1}{a_{k k}}\sum_{i=1}^da_{k
  i}\mu_i-\frac{1}{a_{k k}}\left(\sum_{a_{k i}<0, i\neq k}a_{k
    i}b_i+\sum_{a_{k i}\geq 0, i\neq k}a_{k i}a_i\right)\,.\nonumber 
\end{eqnarray}

\paragraph{Remark}
Since $R_k(\mathcal{B})>0$ for all $k$, $\tau(t)>-\infty$ almost surely for all
$t\in\Z$. This implies that the truncated multivariate normal case falls under
Theorem's \ref{p31} hypothesis and our algorithm works.

\paragraph{Back to one dimension. Minimum of truncated normals}
We prove Proposition~\ref{p22}.  To simplify notation write $g_{a,b}(x;\mu)$
instead of $g_{a,b}(x;\mu,1)$.
Observing that 
\[
g_{a,b}(x;\mu,\sigma^2)=
\frac{1}{\sigma}g_{\frac{a}{\sigma},\frac{b}{\sigma}}
\left(\frac{x}{\sigma},\frac{\mu}{\sigma}\right)  
\]
it is sufficient to prove \reff{Marzo1} for $\sigma=1$.  The proof is based in
the following elementary lemmas.

\begin{lemma}
\label{S14}
 Let $\mu_1<\mu_2$. Then
\begin{eqnarray}
\min(g_{a,b}(x;\mu_1),g_{a,b}(x;\mu_2))&=&
g_{a,b}(x;\mu_1)\one\{x(\mu_1,\mu_2)\leq x\}\nonumber\\
&&\qquad+\;g_{a,b}(x;\mu_2)\one\{x(\mu_1,\mu_2)>x\},
\end{eqnarray}
where
\begin{equation}
x(\mu_1,\mu_2)=
\frac{\mu_1+\mu_2}{2}-\frac{1}{\mu_2-\mu_1}
\ln\left(\frac{A_{a,b}(\mu_1)}{A_{a,b}(\mu_2)}\right).
\end{equation}
\end{lemma}

\begin{lemma}
\label{S16}
For all $\mu\in (\mu^-,\mu^+)$ it holds
\begin{eqnarray}
\label{S17}
\frac{\mu+\mu^+}{2}-\frac{1}{\mu^+-\mu}
\ln\left(\frac{A_{a,b}(\mu)}{A_{a,b}(\mu^+)}\right)
\;>\;\frac{\mu^-+\mu}{2}-\frac{1}{\mu-\mu^-}
\ln\left(\frac{A_{a,b}(\mu^-)}{A_{a,b}(\mu)}\right). 
\end{eqnarray}
\end{lemma}

\paragraph{Proof of Proposition \ref{p22}}
Let $\mu\in (\mu^-,\mu^+)$. By Lemma \ref{S14} we have
\begin{eqnarray}
g_{a,b}(x;\mu^-)\leq g_{a,b}(x;\mu) & \Leftrightarrow &
\frac{\mu^-+\mu}{2}-\frac{1}{\mu-\mu^-}
\ln\left(\frac{A_{a,b}(\mu^-)}{A_{a,b}(\mu)}\right)\leq x,\label{p17}\\
g_{a,b}(x;\mu)\leq g_{a,b}(x;\mu^+) & \Leftrightarrow &
\frac{\mu+\mu^+}{2}-\frac{1}{\mu^+-\mu}
\ln\left(\frac{A_{a,b}(\mu)}{A_{a,b}(\mu^+)}\right)\leq x\,. 
\label{p18}
\end{eqnarray}
Hence $g_{a,b}(x;\mu)\leq\min(g_{a,b}(x;\mu^-),g_{a,b}(x;\mu^+))$ if
and only if
\begin{eqnarray}
\label{S18}
x\in\left[\frac{\mu+\mu^+}{2}-\frac{1}{\mu^+-\mu}
\ln\left(\frac{A_{a,b}(\mu)}{A_{a,b}(\mu^+)}\right),
\frac{\mu^-+\mu}{2}-\frac{1}{\mu-\mu^-}
\ln\left(\frac{A_{a,b}(\mu^-)}{A_{a,b}(\mu)}\right)\right]\,.
\end{eqnarray}
But by Lemma \ref{S16} the interval in \reff{S18} is empty for all $\mu\in
(\mu^-,\mu^+)$. This implies that\hfill\break
$\min(g_{a,b}(x;\mu^-),g_{a,b}(x;\mu^+))<g_{a,b}(x;\mu)$, from where
\[\inf_{\mu\in I}g_{a,b}(x,\mu,\sigma^2)
=\min(g_{a,b}(x;\mu^-,\sigma^2),g_{a,b}(x;\mu^+,\sigma^2))\]
holds for $\sigma=1$.
The corresponding identity (\ref{Marzo1}) follows by applying Lemma \ref{S14} to
$\mu_1=\mu^-$ and $\mu_2=\mu^+$.
\hfill{$\square$}

\section{Comparisons and comments on other methods}
\label{s6}

In this section we compare our method with the rejection method in some
examples. The conclusion is that the rejection method may be better than ours in
dimension 2 for some regions but ours becomes better and better as dimension
increases.  Then we show why the method does not work when the region is not a
box. This discards the following tempting approach: multiply the target vector
by a matrix to obtain a vector with iid coordinates. The transformed vector is
much easier to simulate; the problem is that it is now conditioned to a
transformed region.  When the region is not a box, the conditioned vector is not
an iid vector and a coupling must be performed. We show here that in general the
corresponding coupling event has probability zero.

\paragraph{The rejection method} We compare our algorithm with the following
rejection algorithm: simulate a uniformly distributed vector $(\bx,u)$ in
$\cB\times [0,\max_{\bx\in \cB} g(\bx)]$. If $u<g(\bx)$, then accept $\bx$.  We
use $d+1$ one-dimensional uniform random variables for each attempt of the
rejection algorithm.  Let $N$ be the expected number of one-dimensional uniform
random variables generated by our perfect simulation algorithm, each counted the
number of times that it is used.  Let $M$ be the expected number of
one-dimensional uniform random variables needed by the rejection algorithm to
produce an accepted value. We divide the experiments in two parts: $d=2$ and
$d\ge 2$.

\begin{figure}[ht]
\begin{center}
  \subfigure[(ii) type
  $1$.]{\includegraphics[width=.32\textwidth]{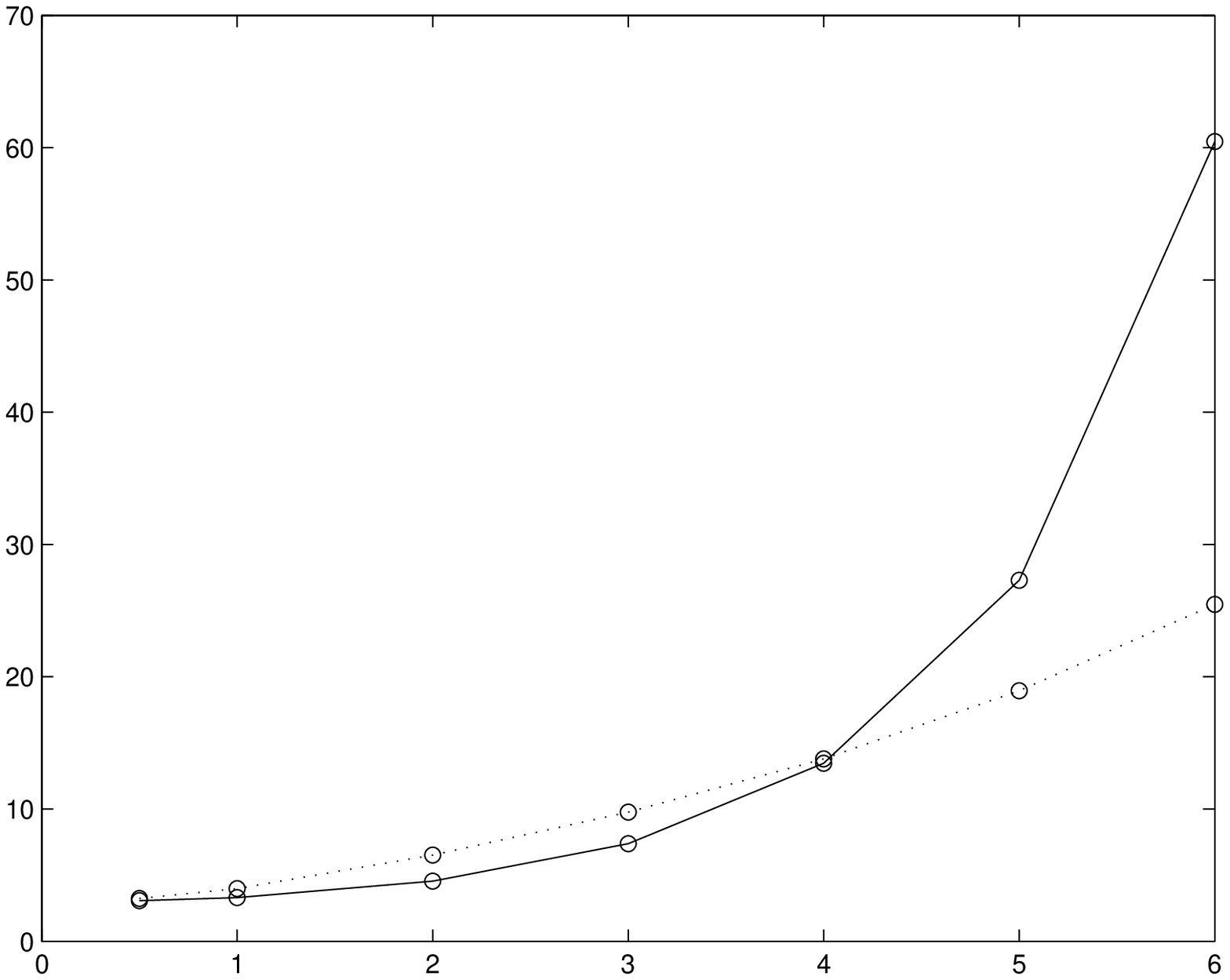}} \subfigure[(ii)
  type $2$.]{\includegraphics[width=.32\textwidth]{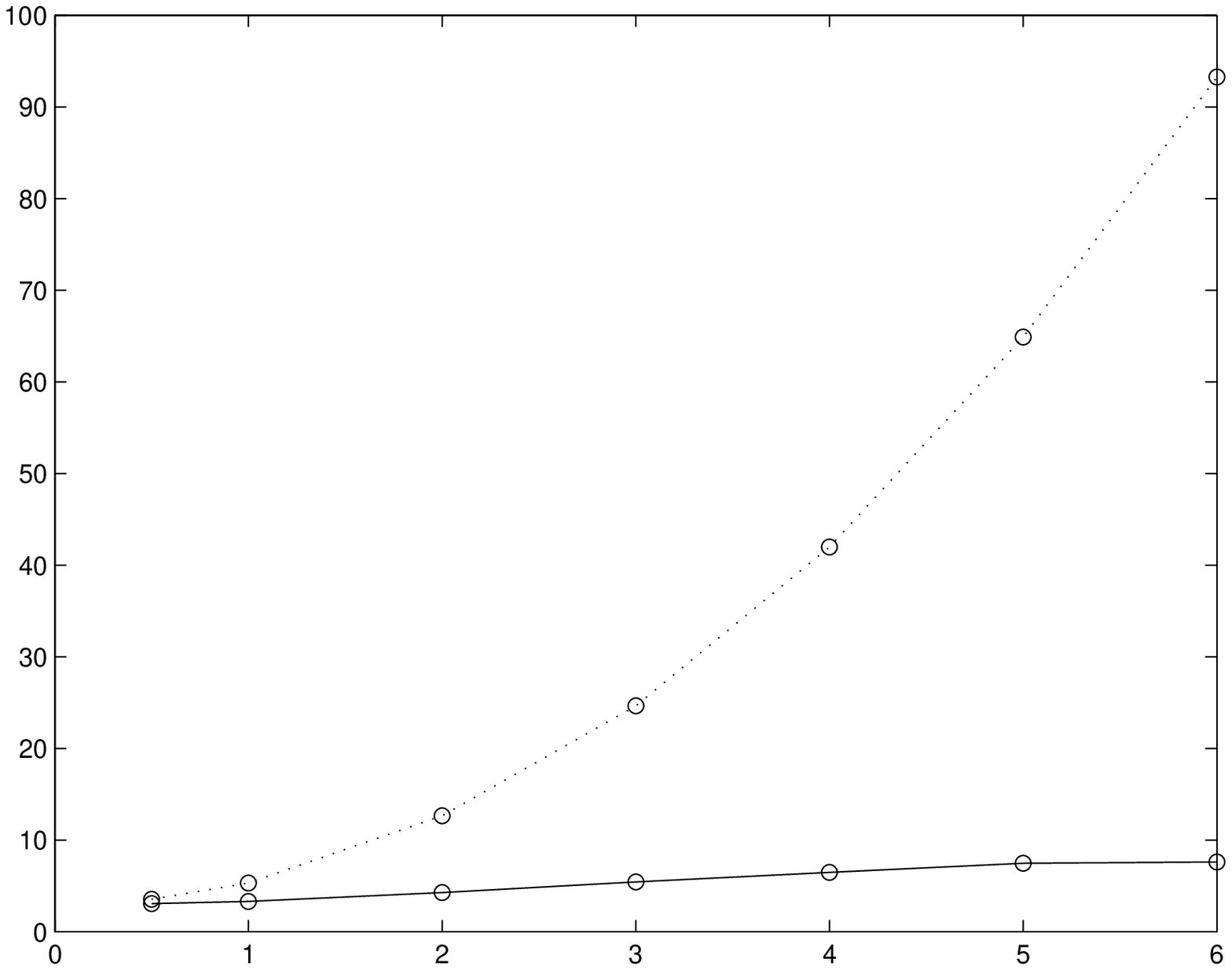}}
\end{center}
\caption{\label{dos} $d=2$ case (ii): $(r,N(r))$ (in solid line) and $(r,M(r))$
  (in dashed line) for boxes of both types. For type $2$ our method is faster
  than the rejection method, but for type $1$ this is true when $r\leq 4$.}
\end{figure}

$[d=2]$.  Let $\mu=(0,0)$ and
$\bSigma=\left(\begin{array}{cc}1&12/5\\12/5&9\end{array}\right)$. We consider
two cases. (i) boxes $\cB(x_1,x_2)=[0,1]\times[0,1]+(x_1,x_2)$, $-4\leq
x_1\leq 4$ and $0\leq x_2\leq 4$. (ii) boxes $\cB(r)=[0,r]\times[0,r]$ (type $1$)
and $\cB(r)=[-r,0]\times[0,r]$ (type $2$), $r=1/2,1,2,3,4,5,6$.  The results are
shown in Table~\ref{uno} for (i) and in Figure \ref{dos} for (ii).

\begin{table}[ht]
\begin{eqnarray*}
\begin{array}{c|ccccccccc}
x_2\setminus x_1&-4&-3&-2&-1&0&1&2&3&4\\
\hline
0&65&37&16&5&4&12&29&55&86\\
\hline
1&76&45&22&8&4&8&22&45&75\\
\hline
2&88&55&28&11&4&6&17&37&64\\
\hline
3&100&65&36&15&5&5&12&30&54\\
\hline
4&112&77&45&21&7&4&9&23&45
\end{array}
\end{eqnarray*}
\caption{\label{uno} $d=2$ case (i): 
  values of $M(x_1,x_2)$, mean number of one-dimensional uniform random
  variables used by the rejection method. The corresponding number for our
  method, $N(x_1,x_2)\in [3.2,3.3]$ is better for  all boxes.} 
\end{table}

\begin{figure}[ht]
\begin{center}
\subfigure[]{\includegraphics[width=.32\textwidth]{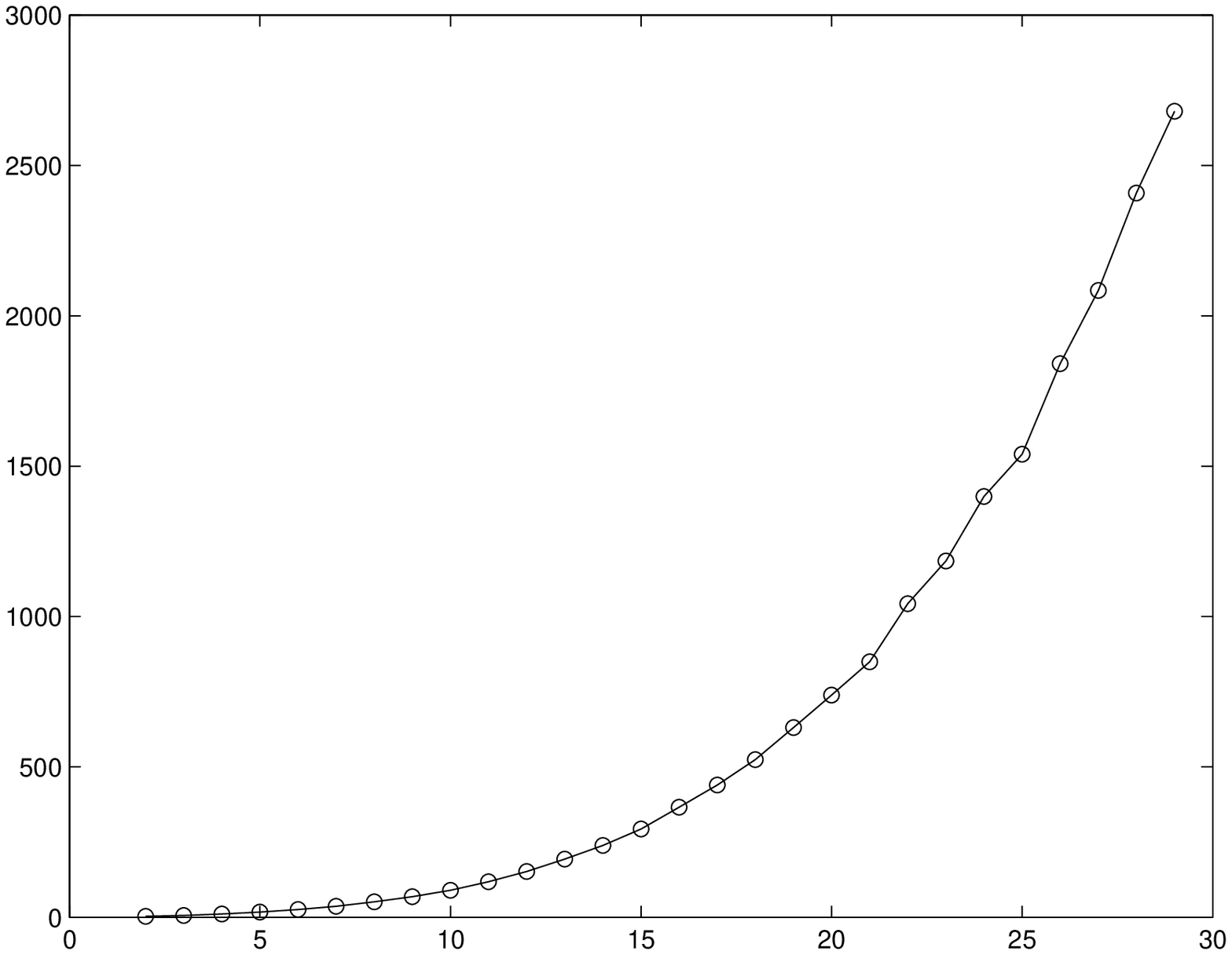}}
\subfigure[]{\includegraphics[width=.32\textwidth]{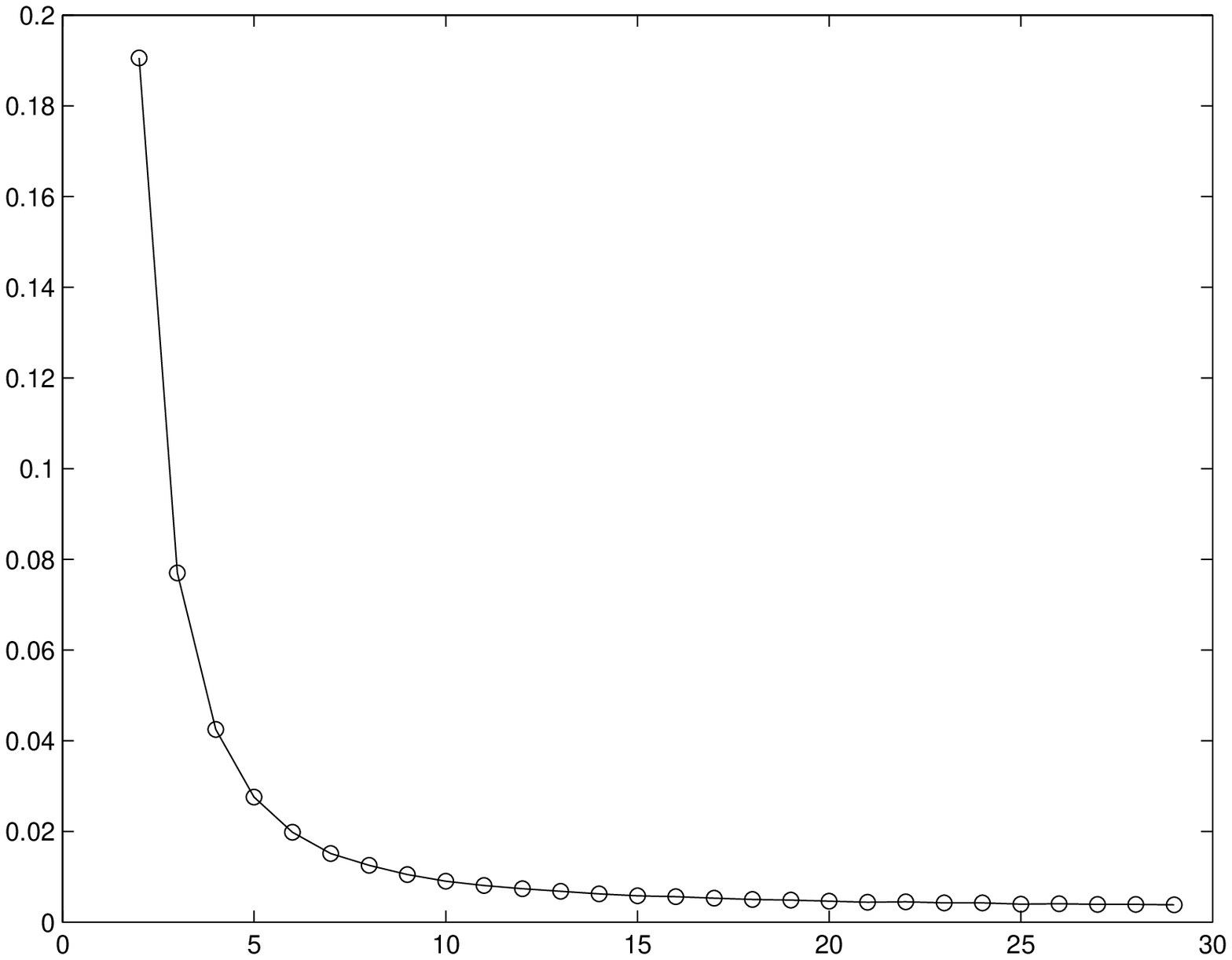}}
\subfigure[]{\includegraphics[width=.32\textwidth]{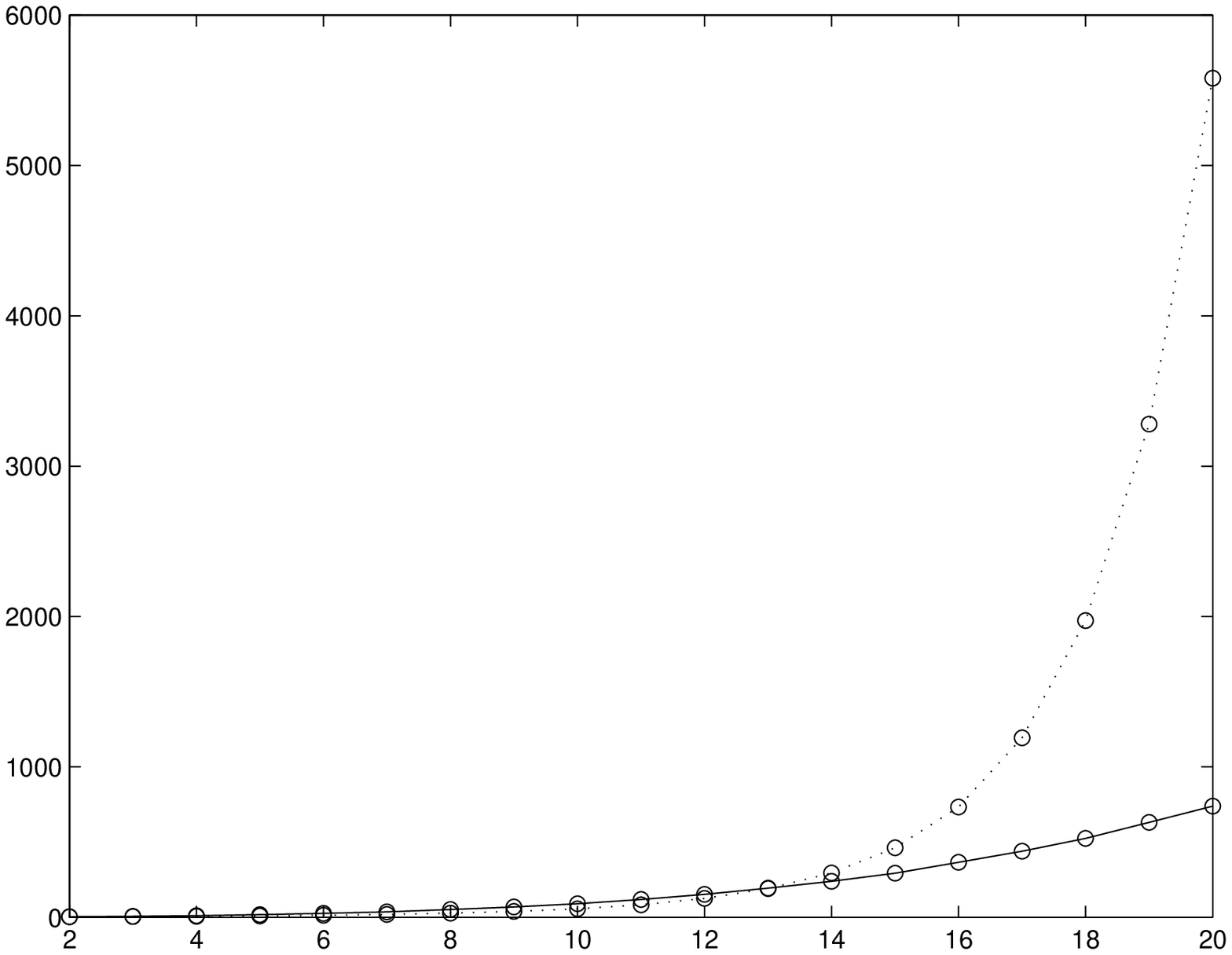}}
\end{center}
\caption{\label{tres} Case (ii) $d\ge 2$. (a) $(d,N(d))$; (b)
  $(d,N(d)/d^4)$; (c) $(d,N(d))$ (solid line) and $(d,M(d))$ (dotted line). It
  seems that $N(d)$ grows as $o(d^4)$ and $M(d)$ as $(d+1)e^{ad}$, for some
  $a>0$. For $d\geq 14$ our method is faster than the rejection method.}
\end{figure}

\clearpage
$[d\ge 2]$ example 1. Fix $\mu=(0,\dots,0)$ and
$\bSigma^{-1}=\frac12I+\frac12\one\one'$, where $\mathbf I$ is the identity matrix and
$\one'=(1,\dots,1)$. We consider two cases: (i) boxes
$\cB=\left[0,\frac12\right]^d$, $d=2,\dots,29$ and (ii) boxes
$\cB=\left[\frac12,1\right]^d$, $d=2,\dots,14$.  The results are shown in Figure
\ref{tres} and Table \ref{cuatro} for (i) and (ii), respectively.

\begin{table}[ht]
\begin{eqnarray*}
\begin{array}{c|ccccccccccccc}
d&2&3&4&5&6&7&8&9&10&11&12&13&14\\
\hline
N&3&6&11&17&27&36&48&60&84&106&147&164&208\\
\hline
M&5&10&22&51&133&372&1141&3810&13736&53356&223520&981460&4858138
\end{array}
\end{eqnarray*}
\caption{\label{cuatro} Case (i) $d\ge 2$. $N$ is the mean number of uniforms used by
  our method and $M$ by the rejection method.} 
\end{table}

\begin{figure}[ht]
\begin{center}
\subfigure[]{\includegraphics[width=.32\textwidth]{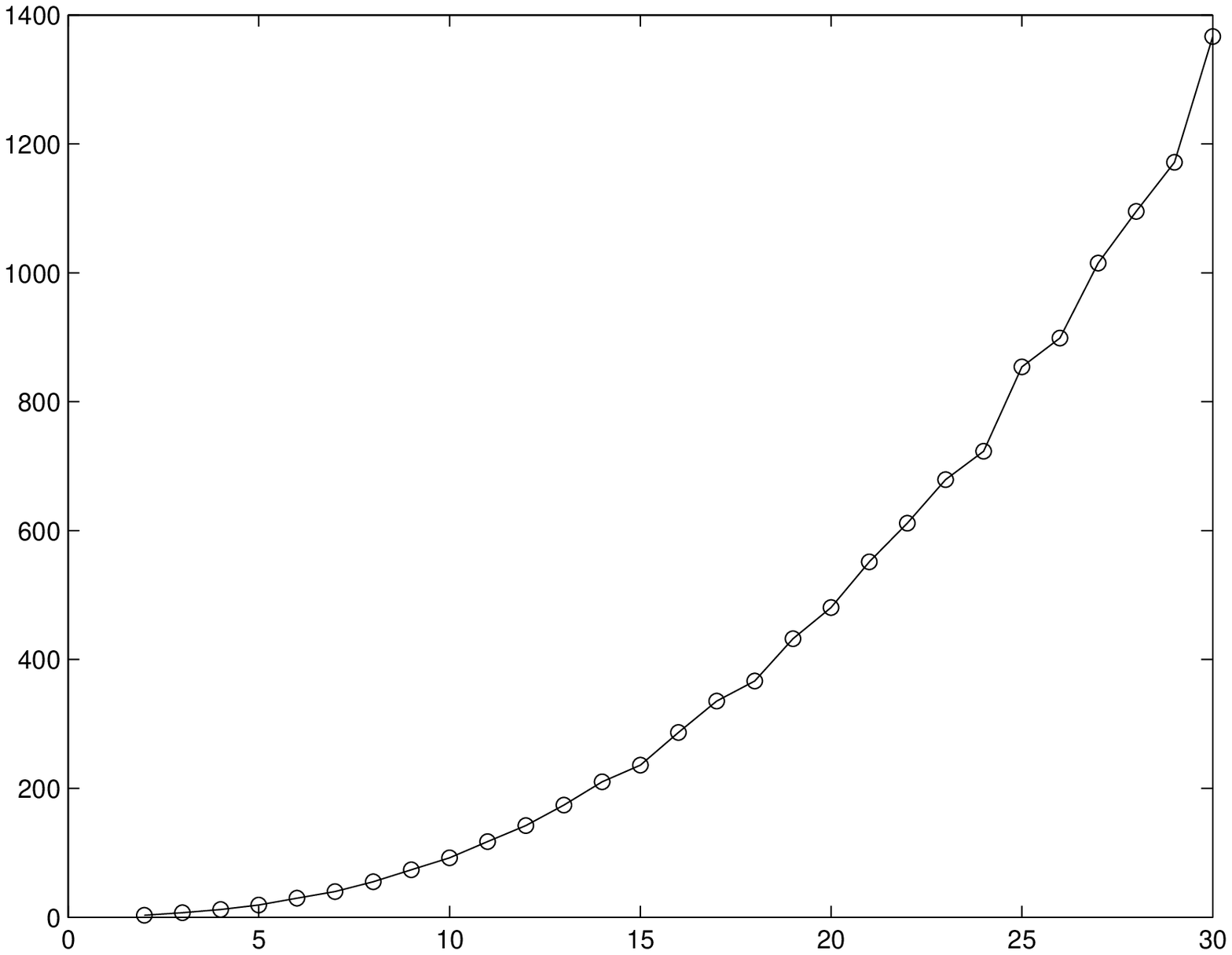}}
\subfigure[]{\includegraphics[width=.32\textwidth]{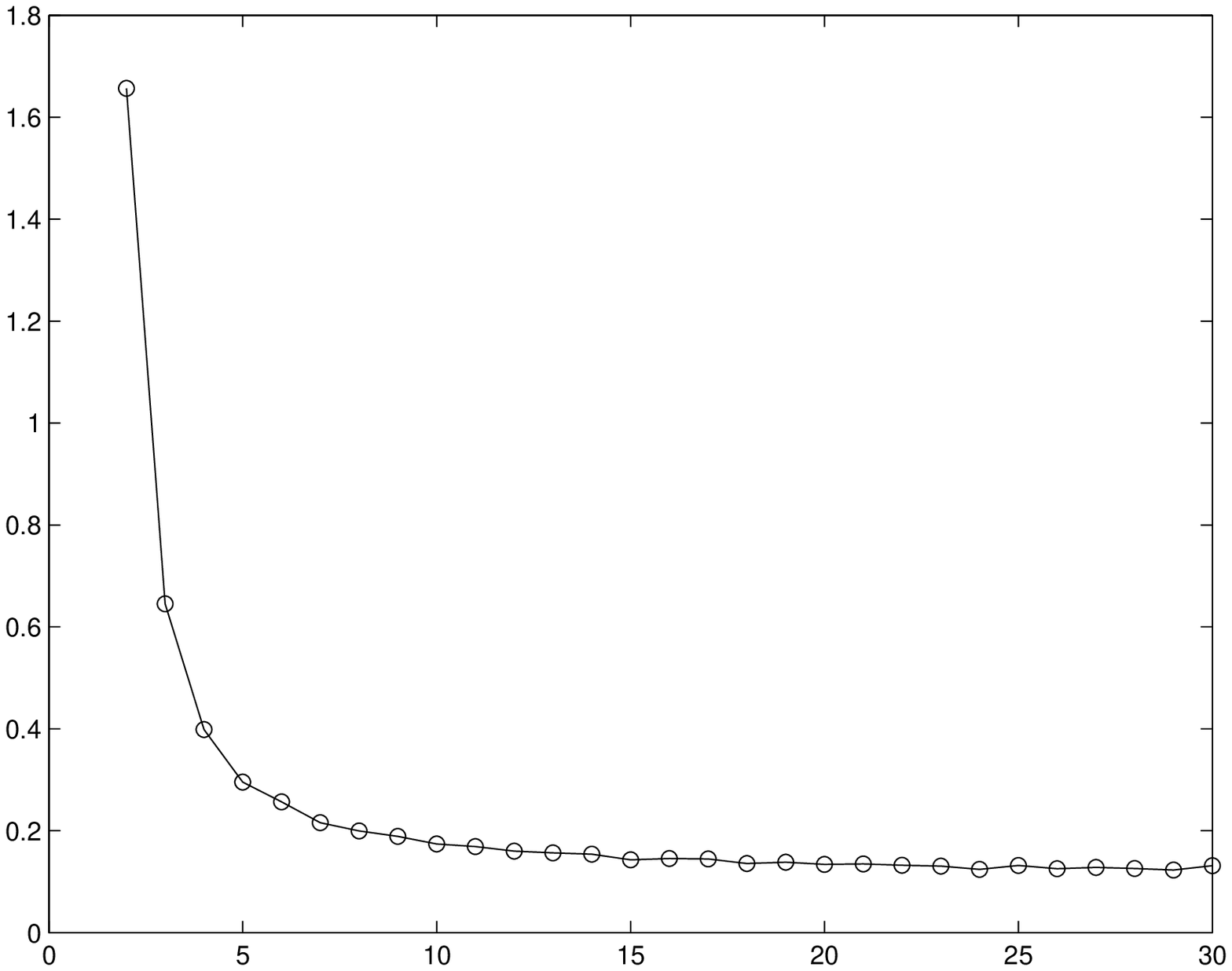}}
\subfigure[]{\includegraphics[width=.32\textwidth]{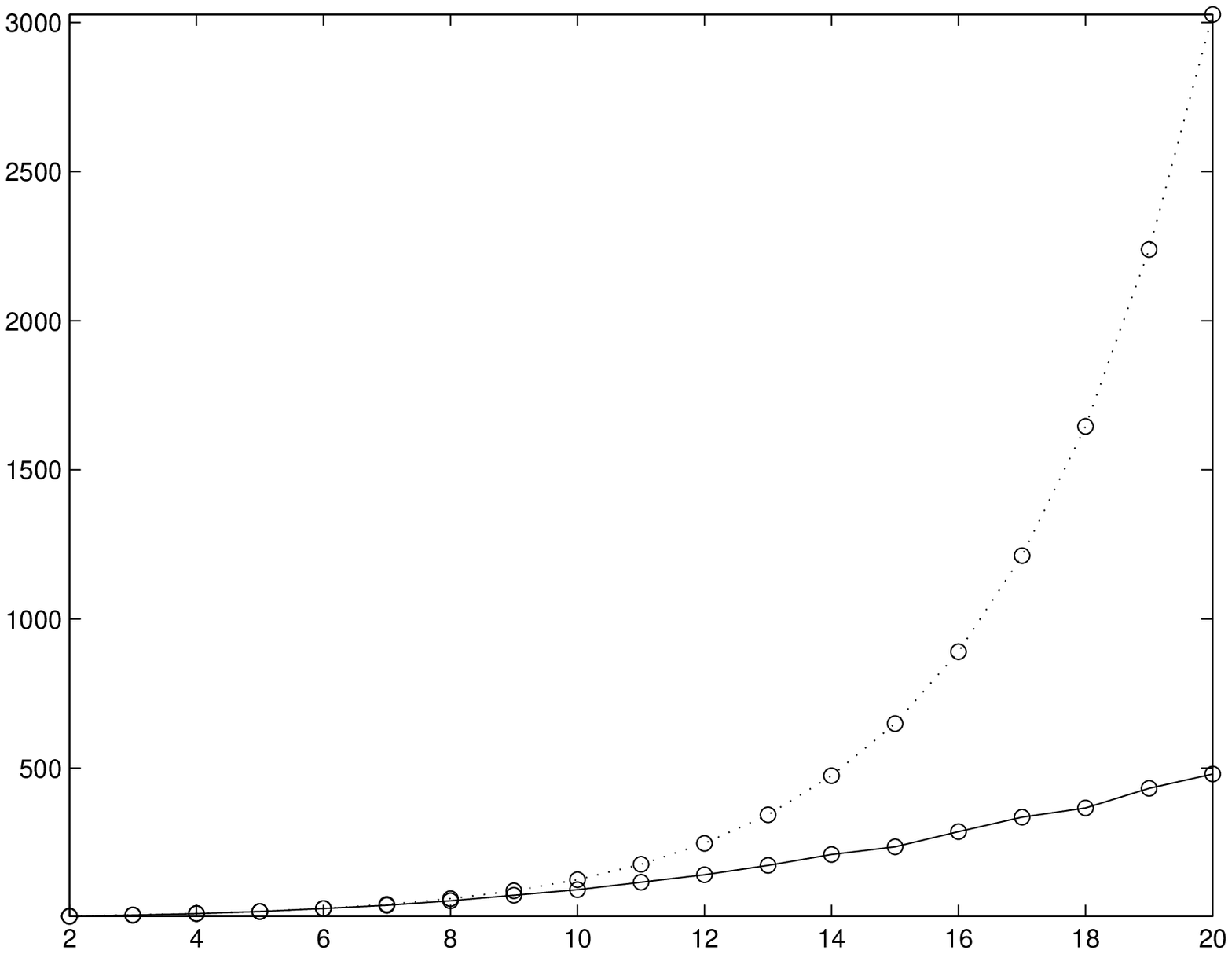}}
\end{center}
\caption{\label{cinco} (a) $(d,N(d))$; (b) $(d,N(d)/(d\log d)^2)$; (c)
  $(d,N(d))$ (in solid line) and $(d,M(d))$ (in dotted line). It seems that
  $N(d)$ grows as $(d \log d)^2$ and $M(d)$ as $(d+1)e^{ad}$, for some
  $a>0$. Our method is faster than the rejection method for all $d$.} 
\end{figure}

$[d\ge 2]$ example $2$.  Let $\mu=(0,\dots,0)$ and $\bSigma^{-1}$ be the matrix,
with ``neighbor structure'', $\bSigma^{-1}_{i,j}=\rho^{|j-i|}\one\{|j-i|\leq
1\}$, with $\rho=\frac12$. We consider boxes of the form $\cB=[0,1]^d$. The
results are shown in Figure \ref{cinco}.

\paragraph{Strongly correlated variables}

The algorithm slows down very fast with the dimension when the variables are
strongly correlated. An Associate Editor proposes to consider a Gaussian vector
truncated to the box $I^d = [0,1]^d$, with mean $\mu=0$ and covariance matrix
\begin{equation}
  \label{b1}
  \bSigma
  =\varepsilon \mathbf I +(1-\varepsilon){\mathbf 1}{\mathbf 1}'
\end{equation}
where $\mathbf I$ is the identity matrix and $\mathbf 1'=(1,\dots,1)$, for small
values of $\epsilon$. As a measure of the speed of the algorithm we have
computed the coefficient $R$, shown in Table~\ref{x2}, as function of $\epsilon$
and the dimension $d$.
\begin{table}[ht]
\begin{eqnarray*}
\begin{array}{c|ccccc}
  \varepsilon\setminus d&2&4&8&16&32\\
  \hline
  0.1&0.5139&0.3446&0.2792&0.2507&0.2375\\
  \hline
  0.01&0.8753\cdot10^{-3}&0.3121\cdot 10^{-4}&0.5969\cdot 10^{-5}&0.2615\cdot 10^{-5}&0.1731\cdot 10^{-5}
\end{array}
\end{eqnarray*}
\caption{\label{x2} Coefficient $R$ for $\mu=0$, $\bSigma
  =\varepsilon \mathbf I +(1-\varepsilon){\mathbf 1}{\mathbf 1}'$ where ${\mathbf I}$ is the
  identity matrix and $\mathbf 1'=(1,\dots,1)$, and the box is $[0,1]^d$.} 
\end{table}
In $d=2$ the coupling time $\tau(0)$ for this example is a geometric with mean
$1/R$, because a it suffices that one of the uniforms is smaller than $R$ to get
the coupling in the next step. In general, as explained in Lemma \ref{s22},
$\tau(0)\le \tau^0(0)$, whose expectation is of the order $(1/R)^{d-1}$.

\paragraph{Other regions}
Let $\bY\sim\cN(0,\bSigma)$ be a $d$-dimensional multinormal random vector with
zero mean and nonsingular covariance matrix $\bSigma$. Let
$\cB=[a_1,b_1]\times\cdots\times[a_d,b_d]$.  If we rotate and scale the
coordinates to obtain iid $\cN(0,1)$'s, i.e., if we put $\bX=\bSigma^{-1/2}\bY$,
then the box constraints are of the form $\ba\leq \bSigma^{1/2}\bx\leq \bb$,
where $\ba=(a_1,\dots,a_d)$ and $\bb=(b_1,\dots,b_d)$.  When the constraints are
of this form, the conditional distribution of $X_k|\bX^k=\bx^k$ has a
probability density function whose support depends on the condition $\bx^k$. In
such case the maximal coupling in our approach may have a coupling event with
zero probability. Indeed, it could happen that 
\[\inf_{\bx} g_{X_k|\bX=\bx}(x_k)=0.\]
This implies that $R_k(x_k)=\int_{-\infty}^{x_k}\inf_{\bx}
g_{X_k|\bX=\bx}(y)dy=0$ for all $x_k$.

To see this in $d=2$, observe that the transformation takes the
rectangular box into a parallelogram $\cF$ with sides not in general
parallel to the coordinate axes. We need to simulate the standard
bi-dimensional normal $\bX$ conditioned to $\cF$. It is still true that
the law of $X_1$ conditioned on $X_2=y$ is normal truncated to the
slice of the parallelogram $\cF^y:=\{(x_1,x_2)\in \cF:\, x_2=y\}$. It
is a matter of geometry that there are different $y$ and $y'$ such that
$\cF^y\cap \cF^{y'}=\emptyset$; this implies the infimum above is zero.

Take for instance $\mu=(0,0)$, $\bSigma=\left(\begin{array}{cc}5&-4\cr
    -4&5\end{array}\right)$ and $\bY$ conditioned to the box
$\cB=[0,1]\times[0,1]$. Then $\bSigma^{-1/2}=\frac13\left(\begin{array}{cc}2&1\cr
    1&2\end{array}\right)$, and the transformed box is a parallelogram $\cF$ with
vertices $(0,0),\, (2/3,1/2),\, (1/3,2/3),\, (1,1)$.

\begin{figure}[ht]
\begin{center}
\input{HR.pstex_t}
\end{center}
\caption{In a parallelogram region the minimal value $\inf_{\bx}
  g_{X_1|\bX=\bx}(x_1)$ is zero. The truncated distribution conditioned to the
  second coordinate being $y$ and $y'$ respectively have disjoint supports
  $\cF^y$ and $\cF^{y'}$, which are the projections on the horizontal axis of
  the green and red segments, respectively.}
\end{figure}
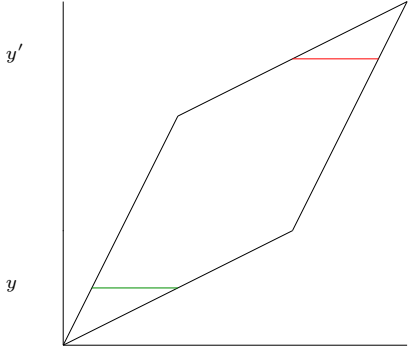

\clearpage

\section{Final remarks}

We have proposed a construction (simulation) of a multidimensional random vector
conditioned to a bounded box using a version of the coupling from the past
algorithm of the Gibbs sampler dynamics. In the case of normal distribution, we
have given a explicit construction taking advantage of the bell shape of the
density and the fact that the distribution of one coordinate given the others is
a truncated normal with variance independent of the other coordinates. This
allows to define the coupling set $C$ in function of the extreme possible values
of the mean and to obtain positive probability for each coupling set $C$ as soon
as the box is bounded. 

One of the consequences of the approach is to show that the Gibbs sampler is
ergodic in the box if the coupling event has probability uniformly bounded below
in each one of the coordinates. This corresponds to the condition $R_k(\cB)>0$
for all coordinate $k$. When the uniform random variable used to update site $k$
at some time falls below $R_k(\cB)>0$, this coordinate coincides for all
processes.

The method can be used to show ergodicity of the process in infinite volume when
the inverse of the covariance matrix has a neighbor structure; for instance for
tridiagonal matrices. As an example, consider the formal density in
$[a,b]^{\Z^d}$:
\[
g(\bx) \sim \exp \Bigl(-\beta \sum_{|i-j|=1} (x_i-x_j)^2 \Bigr) 
\] 
which corresponds to an infinite volume Gaussian field, each coordinate
truncated to $[a,b]\subset \R$; $\beta>0$ is a parameter (inverse temperature
in Statistical Mechanics). The Gibbs sampler can be performed in continuous time
and the coupling can be done as in the finite case. However the updating of site
$k$ depends only on the configurations in sites $j$ such that $|j-k|=1$.  Each
site gets a definite value (and not an interval) when the corresponding uniform
random variable falls below a certain value $R=R(\beta,[a,b])$ (corresponding to
$R_k(\cB)$, but now it is constant in the coordinates). To determine the value
of site $k$ at time $t$ one explores backwards the process; calling $U$ the
random variable used for the last time (say $t'$) before $t$ to update site $k$
there are two cases: $U<R$ or $U>R$. If $U<R$, we do not need to go further
back, the value is determined. If $U>R$, we need to know the values of $2d$
neighbors at time $t'$. Repeating the argument, we construct a ``oriented
percolation'' structure which will eventually finish if $R>1-\frac{1}{2d}$. This
value is obtained by dominating the percolation structure with a branching
process which dies out with probability $R$ and produces $2d$ offsprings with
probability $1-R$. The value of $R$ depends on the length of
the interval $[a,b]$ and the strength of the interaction governed by $\beta$. One
can imagine that for small intervals and small $\beta$ things will work.

\section*{Acknowledgements} We thank Christian Robert and Havard Rue for
discussions.

This paper was partially supported by FAPESP, CNPq, CAPES-SECyT, Facultad de
Ingenier\'ia de la Universidad de Buenos Aires.

\end{document}

%% file: HR.pstex_t
\begin{picture}(0,0)%
\includegraphics{HR.pstex}%
\end{picture}%
\setlength{\unitlength}{2368sp}%
\begingroup\makeatletter\ifx\SetFigFont\undefined%
\gdef\SetFigFont#1#2#3#4#5{%
  \reset@font\fontsize{#1}{#2pt}%
  \fontfamily{#3}\fontseries{#4}\fontshape{#5}%
  \selectfont}%
\fi\endgroup%
\begin{picture}(4212,3624)(601,-3973)
\put(601,-3361){\makebox(0,0)[lb]{\smash{\SetFigFont{7}{8.4}{\rmdefault}{\mddefault}{\updefault}{\color[rgb]{0,0,0}$y$}%
}}}
\put(601,-961){\makebox(0,0)[lb]{\smash{\SetFigFont{7}{8.4}{\rmdefault}{\mddefault}{\updefault}{\color[rgb]{0,0,0}$y'$}%
}}}
\end{picture}

%% file: ffg-jap-revised.bbl
\begin{thebibliography}{Dillo 60}



\bibitem{AGR99} Allenby, G. M. and Rossi, P. (1999) Marketing Models of
  consumer Heterogeneity. {\it Journal of Econometrics, \bf 89}, 57-78.

\bibitem{D86} Devroye, L. (1986) {\it Non-Uniform Random Variate Generation},
  Springer-Verlag, New York.

\bibitem{GS90} Gelfand, A. E. and  Smith, A. F. M. (1990) Sampling-Based
  Approaches to Calculating Marginal Densities, {\it Journal of the American
  Statistical Association, \bf 85}, 398-409.

\bibitem{GG84} Geman, S. and Geman, D. (1984) Stochastic Relaxation, Gibbs
  Distribution and the Bayesian Restoration of Images, {\it IEEE Transactions on
    Pattern Analysis and Machine Intelligence, \bf 6}, 721-741.

\bibitem{G91} Geweke, J. (1991) Efficient simulation from the Multivariate
  Normal and Student t-distribution subject to linear constrains. {\it Computer
  Sciences and Statistics}: Proc. 23d Symp. {\it Interface}, 571-577.

\bibitem{GKR97} Geweke, J., Keane, M. and Runkle, D.(1997) Statistical
  inference in the Multinomial Probit Model. {\it Journal of Econometrics, \bf
  80}, 125-165.

 \bibitem{JKB} Johnson,  N. L., Kotz, S. and   Balakrishnan, N. (1994) {\it
     Continuous univariate distributions}, Volume 1, Wiley, New York. 

\bibitem{m99} Moeller, J. (1999) Perfect simulation of conditionally specified
  models. {\it Journal of the Royal Statistical Society, Ser. B, \bf
    61}(1):251--264.

\bibitem{mg98}  Murdoch, D. J. and Green,  P. J. (1998) Exact sampling from a
  continuous state space. {\it Scandinavian Journal of Statistics \bf
    25}(3):483--502.
        
  
\bibitem{pr03} Philippe, A. Robert, C. (2003) Perfect simulation of
  positive Gaussian distributions {\it Statistics and Computing, \bf 13},
  Number 2 179 - 186

\bibitem{PW96} Propp, J. G. and Wilson, D. B. (1996) Exact
  sampling with coupled Markov chains and applications to statistical mechanics.
  Proceedings of the Seventh International Conference on Random Structures and
  Algorithms (Atlanta, GA, 1995).  {\it Random Structures Algorithms \bf 9}, no.
  1-2, 223--252.

\bibitem{R95} Robert, C. (1995) Simulation of truncated normal random
  variables. {\it Statistics and computing \bf 5}, 121-125.

\bibitem{Ru87} Rubin, D. B. (1987) {\it Multiple imputations for Nonresponse in
  Surveys}. New York. Wiley.
 
\bibitem{TW87} Tanner, M. A. and Wong, W. (1987) The Calculation of
  Posterior Distributions by Data Augmentation, {\it Journal of the American
  Statistical Association, \bf 82}, 528-549.

\bibitem{Ta91} Tanner, M. A. (1991) Tools for Statistical Inference:
  Observed Data and Data Augmentation Methods. {\it Lecture Notes in
    Statistics \bf 67}. New York: Springer Verlag.

\bibitem{Tho} Thorisson, H. (2000) {\it Coupling, Stationarity, and
    Regeneration}, Springer - Verlag, New York.


\bibitem{wil} D.~B. Wilson.  \newblock Annotated bibliography of perfectly
  random sampling with \char77arkov chains.  \newblock In D.~Aldous and
  J.~Propp, editors, {\em Microsurveys in Discrete Probability}, volume~41 of
  {\em DIMACS Series in Discrete Mathematics and Theoretical Computer Science},
  pages 209--220. American Mathematical Society, 1998.  \newblock Updated
  versions can be found at {\tt
    http://dimacs.rutgers.edu/\char126dbwilson/exact}.


\end{thebibliography}
